\documentclass[a4paper,12pt]{amsart}

\usepackage{latexsym, amsfonts, amsmath, amsthm, amssymb,verbatim,amscd}

\let\ds=\displaystyle

\newcommand{\AC}{{AC}}
\newcommand{\BV}{{BV}}
\newcommand{\ACHK}{{AC_{HK}}}
\newcommand{\CTPP}{\hbox{$CT\kern-0.2ex{P}\kern-0.2ex{P}$}}

\newcommand{\HP}{\mathcal{H}}  
\newcommand{\NP}{\mathcal{N}}  

\newcommand{\mN}{\mathbb{N}}
\newcommand{\mR}{\mathbb{R}}
\newcommand{\mC}{\mathbb{C}}
\newcommand{\mD}{\mathbb{D}}
\newcommand{\mT}{\mathbb{T}}

\newcommand{\lam}{\boldsymbol{\lambda}}

\newcommand{\Real}{{\rm Re}}
\newcommand{\Imag}{{\rm Im}}

\newcommand{\set}[1]{\left\{#1\right\}}
\newcommand{\abs}[1]{\left\lvert#1\right\rvert}

\newcommand{\norm}[1]{\left\lVert#1\right\rVert}
\newcommand{\ssnorm}[1]{\lVert #1 \rVert}
\newcommand{\norminf}[1]{\left\lVert#1\right\rVert_\infty}
\newcommand{\normbv}[1]{\left\lVert#1\right\rVert_{\BV(\sigma)}}
\newcommand{\normbvj}[1]{\left\lVert#1\right\rVert_{\BV(J)}}

\def\ipr<#1,#2>{\langle #1,#2\rangle}

\newcommand{\Range}{{\rm Range}}

\newcommand{\Proj}{{\rm Proj}}

\DeclareMathOperator{\var}{var}
\DeclareMathOperator{\cvar}{\rm cvar}
\DeclareMathOperator{\vf}{vf}

\DeclareMathOperator{\supp}{supp}

\newtheorem{thm}{Theorem}[section]
\newtheorem{cor}[thm]{Corollary}
\newtheorem{lem}[thm]{Lemma}
\newtheorem{prop}[thm]{Proposition}
\newtheorem{defn}[thm]{Definition}
\newtheorem{ex}[thm]{Example}

\newtheorem{quest}[thm]{Question}

\newcommand{\booktitle}[1]{\textit{#1}}
\newcommand{\journalname}[1]{\textrm{#1}}

\title{$\AC(\sigma)$ operators}
\author{Brenden Ashton} 
\email{brenden.ashton@silverbrookresearch.com}

\author{Ian Doust}
\address{School of Mathematics and Statistics\\
University of New South Wales\\
UNSW Sydney 2052 Australia}%
\email{i.doust@unsw.edu.au}%
\date{}
\subjclass{47B40}%
\keywords{Functions of bounded variation, absolutely continuous
functions, functional calculus, well-bounded operators,
$AC$-operators.}
\date{}

\begin{document}
%
%
\begin{abstract}
In this paper we present a new extension of the theory of well-bounded operators to cover operators with complex
spectrum. In previous work a new concept of the class of absolutely continuous functions on a nonempty compact
subset $\sigma$ of the plane, denoted $\AC(\sigma)$, was introduced.  An $\AC(\sigma)$ operator is one which
admits a functional calculus for this algebra of functions. The class of $\AC(\sigma)$ operators includes all of
the well-bounded operators and trigonometrically well-bounded operators, as well as all scalar-type spectral
operators, but is strictly smaller than Berkson and Gillespie's class of $\AC$ operators. This paper develops
the spectral properties of $\AC(\sigma)$ operators and surveys some of the problems which remain in extending
results from the theory of well-bounded operators.
\end{abstract}

\maketitle

%
%

\section{Introduction}\label{intro}
A Banach space operator is said to be well-bounded if it admits a functional calculus for $\AC(J)$, the algebra of absolutely continuous functions on some compact interval $J \subseteq \mR$. The motivation for the introduction of this class was to provide a theory which extended the spectral representation results which apply to self-adjoint operators to Banach space operators which may possess a conditionally rather than unconditionally convergent spectral expansion. Smart and Ringrose 
\cite{dS,jR,jR2} showed that well-bounded operators always have an integral representation with respect to a family of projections known as a
decomposition of the identity. The usefulness of this most general form of the theory is somewhat restricted however since the decomposition of the
identity acts on the dual of the underlying Banach space and is in general not unique (see \cite{hD} for examples of this non-uniqueness).

In \cite{BD} a subclass of the well-bounded operators, the 
well-bounded operators of type~(B), were introduced. The type~(B) well-bounded
operators, which include those well-bounded operators acting on
reflexive spaces, possess a theory of integration with respect to a
family of projections which act on the original space. This family
of projections, known as the spectral family, is uniquely determined
by the operator. The integration theory provides an extension of the
$\AC(J)$ functional calculus to a $\BV(J)$ functional calculus where
$\BV(J)$ is the algebra of functions of bounded variation on the
interval $J$.

As is the case for a self-adjoint operator, the spectrum of a well-bounded operator must lie in the real line.
The main obstacle to overcome if one wishes to extend the theory
of well-bounded operators to cover operators whose spectrum may
not lie in the real line, is that of obtaining a suitable
concept of bounded variation for functions defined on a subset
of the plane. Many such concepts exist in the literature. In
\cite{BG}, Berkson and Gillespie used a notion of variation
ascribed to Hardy and Krause to define the $\AC$ operators.
These are the operators which have an $\AC_{HK}(J \times K)$
functional calculus where $\AC_{HK}(J \times K)$ is the algebra
of absolutely continuous functions in the sense of Hardy and
Krause defined on a rectangle $J \times K \subset \mR^2 \cong
\mC$. They showed \cite[Theorem 5]{BG} that an operator $T \in
B(X)$ is an $\AC$ operator if and only if $T = R + i S$ where
$R$ and $S$ are commuting well-bounded operators. In \cite{BDG}
it is shown that this splitting is not necessarily unique.
Furthermore even if $T$ is an $\AC$ operator on a Hilbert space
$H$, it does not necessarily follow that $\alpha T$ is an $\AC$
operator for all $\alpha \in \mathbb{C}$. On the positive side,
the $\AC$ operators include the trigonometrically well-bounded
operators which have found important applications in harmonic analysis and
differential equations (see \cite{BG2} and \cite{BG3}). An
operator $T \in B(X)$ is said to be trigonometrically
well-bounded if there exists a type~(B) well-bounded operator $A
\in B(X)$ such that 
$T = \exp(i A)$. 

One of the problems in the theory well-bounded and $\AC$ operators is that
the functional calculus of these operators is based on an algebra of
functions whose domain is either an interval in the real axis or a rectangle
in the plane. From an operator theory point of view, a much more natural
domain is the spectrum, or at least a neighbourhood of the spectrum.
Secondly, as we have already mentioned, the class of $\AC$ operators is not
closed under multiplication by scalars. This is also undesirable, 
since if one has structural information about an
operator $T$, this clearly gives similar information about $\alpha T$. To
overcome these problems, in \cite{AD} we defined $\AC(\sigma)$, the
Banach algebra of absolutely continuous functions whose domain 
is some compact set $\sigma$ in
the plane. In this paper we look at those operators which have an
$\AC(\sigma)$ functional calculus, which we call $\AC(\sigma)$ operators.

Section~2 summarizes some of the main results from \cite{AD}
concerning the function algebras $\BV(\sigma)$ and
$\AC(\sigma)$. The question as to how one may patch together absolutely 
continuous functions defined on different domains is addressed in Section~3. 
These results will be needed in order to show that $\AC(\sigma)$ operators are 
decomposable in the sense of \cite{CF2}.

In Section~4 we give some results which illustrate the extent of the
class of $\AC(\sigma)$ operators. In particular, we note that this
class contains all scalar-type spectral operators, all well-bounded
operators and all trigonometrically well-bounded operators.

In Section~\ref{lbl:309} we develop some of the main spectral
properties of $\AC(\sigma)$ operators. Here we show that the
$\AC(\sigma)$ operators form a proper subclass of the $\AC$
operators and hence such operators have a splitting into real and
imaginary well-bounded parts. The natural conjecture that every
$\AC(\sigma)$ operator is in fact an $\AC(\sigma(T))$ operator
remains open. Resolving this question depends on being able to
answer some difficult questions about the relationships between
$\AC(\sigma_1)$ and $\AC(\sigma_2)$ for different compact sets
$\sigma_1$ and $\sigma_2$. These issues are discussed in
Section~\ref{support}.

In Section~\ref{lbl:SpecRes} we examine the case where the
$\AC(\sigma)$ functional calculus for $T$ is weakly compact. In this
case one can construct a family of spectral projections associated
with $T$ which is rich enough to recover $T$ via an integration
process. This `half-plane spectral family' is a generalization of
the spectral family associated with a well-bounded operator of
type~(B). A full integration theory for this class of operators is,
however, yet to be developed. In particular, it is not known whether
one can always extend a weakly compact $\AC(\sigma)$ functional
calculus to a $\BV(\sigma)$ functional calculus. The final section
discusses some of the progress that has been obtained in pursuing
such a theory, and lists some of the major obstacles that remain.

Throughout this paper let $\sigma \subset \mC$ be compact and
non-empty. For a Banach space $X$ we shall denote the bounded linear
operators on $X$ by $B(X)$ and the bounded linear projections on $X$
by $\Proj(X)$. Given $T \in B(X)$ with the single valued extension
property (see \cite{nD}) and $x \in X$ we denote the local spectrum
of $x$ (for $T$) by $\sigma_T(x)$. We shall write $\lam$ for the
identity function $\lam : \sigma \rightarrow \mC, z \mapsto z$.

We would like to thank the referee for their careful reading of the manuscript.

%
%

\section{$\BV(\sigma)$ and $\AC(\sigma)$}\label{bv-definitions}

We shall briefly look at $\BV(\sigma)$ and $\AC(\sigma)$. In
particular we look at how two dimensional variation is defined. More
details may be found in \cite{AD}.

To define two dimensional variation we first need to look at
variation along curves. Let $\Gamma = C([0, 1], \mC)$ be the set of
curves in the plane. Let $\Gamma_L \subset \Gamma$ be the curves
which are piecewise line segments. Let $S = \set{z_i}_{i=1}^n
\subset \mC$. We write $\Pi(S) \in \Gamma_L$ for the (uniform speed)
curve consisting of line segments joining the vertices at $z_1, z_2,
\dots, z_n$ (in the given order). For $\gamma \in \Gamma$ we say
that $\set{s_i}_{i=1}^n \subset \sigma$ is a \emph{partition of
$\gamma$ over $\sigma$} if there exists a partition
$\set{t_i}_{i=1}^n$ of $[0, 1]$ such that $t_1 \leq t_2 \leq \dots
\leq t_n$ and such that $s_i = \gamma(t_i)$ for all $i$. We shall
denote the partitions of $\gamma$ over $\sigma$ by $\Lambda(\gamma,
\sigma)$. For $\gamma \in \Gamma$ and $S \in \Lambda(\gamma,
\sigma)$ we denote by $\gamma_S$ the curve $\Pi(S) \in \Gamma_L$.
The variation along $\gamma \in \Gamma$ for a function $f : \sigma
\rightarrow \mC$ is defined as
\begin{equation} \label{lbl:298}
    \cvar(f, \gamma) = \sup_{\set{s_i}_{i=1}^n \in \Lambda(\gamma,
        \sigma)} \sum_{i=1}^{n-1} \abs{f(s_{i+1}) - f(s_i)}.
\end{equation}

To each curve $\gamma \in \Gamma$ we define a weight factor $\rho$.
For $\gamma \in \Gamma$ and a line $l$ we let $\vf(\gamma, l)$
denote the number of times that $\gamma$ crosses $l$ (for a precise
definition of a crossing see Section~3.1 of \cite{AD}). Set
$\vf(\gamma)$ to be the supremum of $\vf(\gamma, l)$ over all lines
$l$. We set $\rho(\gamma) = \frac{1}{\vf(\gamma)}$. Here we take the
convention that if $\vf(\gamma) = \infty$ then $\rho(\gamma) = 0$. We
can extend the definition of $\rho$ to include functions in $C[a,
b]$ in the obvious way.

The two dimensional variation of a function $f : \sigma
\rightarrow \mC$ is defined to be
\begin{equation} \label{lbl:994}
    \var(f, \sigma) = \sup_{\gamma \in \Gamma}
        \rho(\gamma) \cvar(f, \gamma).
\end{equation}
We have the following properties of two dimensional variation which
were shown in \cite{AD}.

\begin{prop}\label{Gamma-L}
Let $\sigma \subseteq \mC$ be compact, and suppose that $f: \sigma
\to \mC$. Then
 \begin{align*}
    \var(f,\sigma)
        &= \sup_{\gamma \in \Gamma_L} \rho(\gamma) \cvar(f, \gamma)
        \\
        & = \sup \Bigl\{ \rho(\gamma_S) \sum_{i=1}^{n-1} \abs{f(s_{i+1}) - f(s_i)}
                 \,:\, \hbox{$S =  \set{s_i}_{i=1}^n \subseteq \sigma$}
                  \Bigr\}.
 \end{align*}
 \end{prop}

\begin{prop} \label{lbl:541}
Let $\sigma_1 \subset \sigma \subset \mC$ both be compact. Let $f,
g : \sigma \rightarrow \mC$, $k \in \mC$. Then
\begin{enumerate}
    \item $\var(f + g, \sigma) \leq \var(f, \sigma) + \var(g, \sigma)$,
    \item\label{item2} $\var(f g, \sigma) \leq \norminf{f} \var(g, \sigma)
        + \norminf{g} \var(f, \sigma)$,
    \item $\var(k f, \sigma) = \abs{k} \var(f, \sigma)$,
    \item $\var(f, \sigma_1) \leq \var(f, \sigma)$.
\end{enumerate}
\end{prop}

For $f : \sigma \rightarrow \mC$ set
\begin{equation} \label{lbl:268}
    \normbv{f} = \norminf{f} + \var(f, \sigma).
\end{equation}
The functions of bounded variation with domain $\sigma$ are
defined to be
\begin{equation*} \label{lbl:415}
    \BV(\sigma) = \set{f : \sigma \mapsto \mC : \normbv{f} < \infty}.
\end{equation*}


To aid the reader we list here some of the main results from
\cite{AD} and \cite{AD2}. The affine invariance of these algebras
(Theorem \ref{lbl:847} and Proposition \ref{lbl:408}) is one of the
main features of this theory and will be used regularly without
comment.

\begin{prop} \label{lbl:199}
If $\sigma = [a, b]$ is an interval then the above definition of
variation agrees with the usual definition of variation. Hence the
above definition of $\BV(\sigma)$ agrees with the usual definition
of $\BV[a, b]$ when $\sigma = [a, b]$.
\end{prop}

\begin{thm} \label{lbl:333}
Let $\sigma \subset \mC$ be compact. Then $\BV(\sigma)$ is a Banach
algebra using the norm given in Equation \eqref{lbl:268}.
\end{thm}

\begin{thm} \label{lbl:847}
Let $\alpha, \beta \in \mC$ and suppose that $\alpha \neq 0$. Then
$\BV(\sigma) \cong \BV( \alpha \sigma + \beta)$.
\end{thm}

\begin{lem}
Let $f : \sigma \rightarrow \mC$ be a Lipschitz function with
Lipschitz constant $L(f) = \sup_{z, w \in \sigma} \abs{\frac{f(z) -
f(w)}{z - w}}$. Then $\var(f, \sigma) \leq L(f) \var(\lam, \sigma)$.
Hence $f \in \BV(\sigma)$.
\end{lem}

We define $\AC(\sigma)$ as being the subalgebra $\BV(\sigma)$
generated by the functions $1$, $\lam$ and $\overline{\lam}$. (Note
that $\lam$ and $\overline{\lam}$ are always in $\BV(\sigma)$.) We
call functions in $\AC(\sigma)$ the \emph{absolutely continuous
functions with respect to $\sigma$}. By Proposition \ref{lbl:199}
this coincides with the usual notion of absolute continuity if
$\sigma = [a, b] \subset \mR$ is an interval. In \cite{AD} the
following properties of $\AC(\sigma)$ are shown.

\begin{prop} \label{lbl:996}
Let $\sigma = [a, b]$ be a compact interval. Let $g \in
\BV(\sigma) \cap C(\sigma)$. Suppose that $\rho(g) > 0$. Then
$\normbv{f \circ g} \leq \frac{1}{\rho(g)}
\norm{f}_{\BV(g(\sigma))}$ for all $f \in \BV(g(\sigma))$.
\end{prop}

\begin{prop} \label{lbl:408}
Let $\alpha, \beta \in \mC$ and suppose that $\alpha \neq 0$. Then
$\AC(\sigma) \cong \AC( \alpha \sigma + \beta)$.
\end{prop}

\begin{prop} \label{lem:ac compact:4590}
If $f \in \AC(\sigma)$ and $f(z) \ne 0$ on $\sigma$ then
$\frac{1}{f} \in \AC(\sigma)$. Indeed, if $M = \ds \inf_{z \in \sigma} |f(z)|$, then 
$\norm{1/f}_{\AC(\sigma)} \le \ds \frac{1}{M} + \frac{\var(f,\sigma)}{M^2}$.
\end{prop}

We shall also need some properties of $\AC(\sigma)$ and
$\BV(\sigma)$ which were not included in \cite{AD}.

\begin{prop} $\BV(\sigma)$ is a lattice. If $f, g \in \BV(\sigma)$, then
  \[ \hbox{$\normbv{f \lor g} \le \normbv{f} + \normbv{g}$ and
$\normbv{f \land g}  \le \normbv{f} + \normbv{g}$.} \]
\end{prop}

\begin{proof}
Suppose that $\gamma \in \Gamma$ and that $\{s_i\}_{i=1}^n \in
\Lambda(\gamma,\sigma)$. Note that for any $a,a',b,b'$,
  \begin{equation}\label{max-inequality}
   |(a \lor a') - (b \lor b')| 
   \le |(a \lor b) - (a' \lor b)| + |(a' \lor b) - (a' \lor b')| 
   \le |a-a'| + |b-b'|
   \end{equation}
and so
  \[
  \sum_{i=1}^{n-1} |(f \lor g)(s_{i+1}) - (f \lor g)(s_i)|
    \le \sum_{i=1}^{n-1} |f(s_{i+1}) - f(s_i)| + |g(s_{i+1}) - g(s_i)|. \]
Thus
  $ \cvar(f \lor g,\gamma) \le \cvar(f ,\gamma) + \cvar(g,\gamma)
  $
and so
  \begin{align*}
  \normbv{f \lor g}
    &= \norm{f \lor g}_\infty + \sup_\gamma \cvar(f \lor g,\gamma)
    \\
    &\le \norm{f}_\infty + \norm{g}_\infty
       + \sup_\gamma \{ \cvar(f ,\gamma) + \cvar(g,\gamma) \} \\
    &\le \norm{f}_\infty + \sup_\gamma \cvar(f,\gamma) +
          \norm{g}_\infty + \sup_\gamma \cvar(g,\gamma)\\
    &= \normbv{f} + \normbv{g}.
  \end{align*}
The proof for $f \land g$ is almost identical.
\end{proof}

Note that $\BV(\sigma)$ is not a \emph{Banach} lattice, even in the
case $\sigma = [0,1]$.

The set $\CTPP(\sigma)$ of functions on $\sigma$ which are
continuous and piecewise triangularly planar relative to $\sigma$
was introduced in \cite{AD}. It is easy to see that $\CTPP(\sigma)$
is a sublattice of $\BV(\sigma)$.

\begin{cor} $\AC(\sigma)$ is a sublattice of $\BV(\sigma)$.
\end{cor}

\begin{proof}
It suffices to show that if $f,g \in \AC(\sigma)$, then $f \lor g
\in \AC(\sigma)$. Suppose then that $f,g \in \AC(\sigma)$. Then
there exist sequences $\{f_n\},\{g_n\} \subseteq \CTPP(\sigma)$ such
that $f_n \to f$ and $g_n \to g$ in $\BV(\sigma)$. As
$\CTPP(\sigma)$ is a lattice, $f_n \lor g_n \in \CTPP(\sigma)$ for
each $n$ and, using (\ref{max-inequality}), one can see that $(f_n
\lor g_n) \to (f \lor g)$. This implies that $f \lor g$ lies in the
closure of $\CTPP(\sigma)$, namely $\AC(\sigma)$.
\end{proof}

If one wishes to apply the results of local spectral theory, it is
important that $\AC(\sigma)$ forms an admissible algebra of
functions in the sense of Colojoar{\v a} and Foia{\c s} \cite{CF2}.
The first step is to show that $\AC(\sigma)$ admits partitions of
unity.

\begin{lem} $\sigma \subset \mC$ be compact. Then $\AC(\sigma)$ is a
normal algebra. That is, given any finite open cover
$\{U_i\}_{i=1}^n$ of $\sigma$, there exist functions
$\{f_i\}_{i=1}^n \subseteq \AC(\sigma)$ such that
\begin{enumerate}
\item $f_i(\sigma) \subset [0,1]$, for all $1 \le i \le n$,
\item $\mathrm{supp} f_i \subseteq U_i$ for all $1 \le i \le n$,
\item $\sum_{i=1}^n f_i  = 1$ on $\sigma$.
\end{enumerate}
\end{lem}

\begin{proof} This follows from the fact that $C^\infty(\sigma)
\subseteq \AC(\sigma)$ \cite[Proposition 4.7]{AD}. More precisely, let $\{U_i\}_{i=1}^n$ be a finite open cover
of $\sigma$ and let $U = \cup_{i=1}^n U_i$. Choose an open set $V$ with $\sigma \subseteq V \subseteq
\overline{V} \subseteq U$. Then there exist non-negative $f_1,\dots,f_n \in C^\infty(V)$ such that $\sum_{i=1}^n
f_i  = 1$ on $V$ (and hence on $\sigma)$, and $\mathrm{supp} f_i \subseteq U_i$ for all $1 \le i \le n$ (see
\cite[page 44]{LM}).
\end{proof}

For $f \in \AC(\sigma)$ and $\xi \not\in \mathrm{supp} f$, define
  \[ f_\xi(z) = \begin{cases}
         \frac{f(z)}{z-\xi}, &  z \in \sigma\setminus \{\xi\}, \\
         0,                  & z \in \sigma \cap \{\xi\}.
         \end{cases}
         \]
Recall that an algebra $\mathcal{A}$ of functions (defined on some
subset of $\mC$)  is admissible if it contains the polynomials, is
normal, and $f_\xi \in \mathcal{A}$ for all $f \in \mathcal{A}$ and
all $\xi \not\in \mathrm{supp} f$.

\begin{prop} \label{lbl:656}
Let $\sigma \subset \mC$ be compact. Then $\AC(\sigma)$ is an
admissible inverse-closed algebra.
\end{prop}

\begin{proof} All that remains is to show that the last property
hold in $\AC(\sigma)$.  Suppose then that $f \in \AC(\sigma)$ and
$\xi \not\in \mathrm{supp} f$. Given that $\mathrm{supp} f$ is
compact, there exists  $h \in C^\infty(\mC)$ such that $h(z) =
(z-\xi)^{-1}$ on $\mathrm{supp} f$ and $h(z) \equiv 0$ on some
neighbourhood of $\xi$. Again using \cite[Proposition 4.7]{AD} we
have that $h|\sigma \in \AC(\sigma)$ and hence that $f_\xi = f h \in
\AC(\sigma)$.
\end{proof}

%
%

\section{Patching theorems}

The relationship between $\var(f,\sigma_1)$, $\var(f,\sigma_2)$ and
$\var(f,\sigma_1 \cup \sigma_2)$ is in general rather complicated.
The following theorem will allow us to patch together functions
defined on different sets.

\begin{thm}\label{variation-join}
 Suppose that $\sigma_1, \sigma_2 \subseteq \mC$ are
nonempty compact sets which are disjoint except at their boundaries. Suppose
that  $\sigma = \sigma_1 \cup \sigma_2$ is convex. If $f:\sigma \to \mC$,
then
  \[ \max\{\var(f,\sigma_1),\var(f,\sigma_2)\}
     \le \var(f,\sigma)
     \le \var(f,\sigma_1) + \var(f,\sigma_2) \]
and hence
  \[ \max\{\norm{f}_{\BV(\sigma_1)},\norm{f}_{\BV(\sigma_2)} \}
  \le  \normbv{f}
  \le \norm{f}_{\BV(\sigma_1)} + \norm{f}_{\BV(\sigma_2)}. \]
Thus, if $f|\sigma_1 \in \BV(\sigma_1)$ and $f|\sigma_2 \in
\BV(\sigma_2)$, then $f \in \BV(\sigma)$.
\end{thm}

\begin{proof} The left-hand inequalities are obvious.

Note that given any points $z \in \sigma_1 \setminus \sigma_2$ and
$w \in \sigma_2 \setminus \sigma_1$ there exists a point $u$ on the
line joining $z$ and $w$ with $u$ in $\sigma_1 \cap \sigma_2$. To
see this, let $\alpha(t) = (1-t)z + t w$ and let $t_0 = \inf\{ t \in
[0,1] \,:\, \alpha(t) \in \sigma_2\}$. By the convexity of $\sigma$,
$\alpha(t) \in \sigma_1$ for all $0 \le t < t_0$. The closedness of
the subsets then implies that $u = \alpha(t_0)  \in \sigma_1 \cap
\sigma_2$.

Suppose then that $S = \{z_0,z_1,\dots,z_n\} \subseteq \sigma$. For
any $j$ for which $z_j$ and $z_{j+1}$ lie in different subsets, then
using the above remark, expand $S$ to add an extra vertex on the
line joining $z_j$ and $z_{j+1}$ which lies in both $\sigma_1$ and
$\sigma_2$. (Note that the addition of these extra vertices does not
change the value of $\rho(\gamma_S)$ and can only increase the
variation of $f$ between the vertices.) Write the vertices of
$\gamma$ which lie in $\sigma_1$ as $S_1 =
\{z_0^1,z_1^1,\dots,z_{k_1}^1\}$ and those which lie in $\sigma_2$
as $S_2=\{z_0^2,z_1^2,\dots,z_{k_2}^2\}$, preserving the original
ordering. Note that for every $j$, $\{z_j,z_{j+1}\}$ is subset of at
least one of the sets $S_1$ and $S_2$. Thus
  \[ \sum_{j=1}^n |f(z_j)-f(z_{j-1})| \le
     \sum_{i=1}^2 \sum_{j=1}^{k_i} |f(z_j^i)-f(z_{j-1}^i)|
              \]
where an empty sum is interpreted as having value $0$.  Recall that
if $S' \subseteq S$ then $\rho(\gamma_{S'}) \ge \rho(\gamma_S)$.
Thus
  \begin{align*}
      \rho(\gamma_S) \sum_{j=1}^n |f(z_j)-f(z_{j-1})|
      & \le \sum_{i=1}^2 \rho(\gamma_{S_i})
                \sum_{j=1}^{k_i} |f(z_j^i)-f(z_{j-1}^i)| \\
      & \le \sum_{i=1}^2 \rho(\gamma_{S_i}) \cvar(f,\sigma_i) \\
      & \le \sum_{i=1}^2 \var(f,\sigma_i).
  \end{align*}
The results follows on taking a supremum over finite $S \subseteq
\sigma$.
\end{proof}

Note that the convexity of $\sigma$ is vital in
Theorem~\ref{variation-join}. Without this condition it is easy to construct
examples where $\var(f,\sigma_1) + \var(f,\sigma_2) = 0$
for a non constant function $f$.

Later, we will need to show that we can patch two absolutely
continuous functions together. For notational simplicity, the following lemma 
is stated in terms of specific sets $\sigma_1$ and $\sigma_2$, but the affine invariance result (Proposition~\ref{lbl:408}) implies that this immediately also applies to any two rectangles that meet along an edge.

\begin{lem}\label{pasting-lemma}
Suppose that $\sigma_1 = [0,1] \times [0,1]$, that $\sigma_2 = [1,2]
\times [0,1]$ and that $\sigma = \sigma_1 \cup \sigma_2$. Suppose
that $f: \sigma \to \mC$ and that $f_i = f|\sigma_i$ ($i=1,2$). If
$f_1 \in \AC(\sigma_1)$ and $f_2 \in \AC(\sigma_2)$, then $f \in
\AC(\sigma)$ and
  \[ \normbv{f}\le \norm{f_1}_{\BV(\sigma_1)} +
\norm{f_2}_{\BV(\sigma_2)}. \]
\end{lem}

\begin{proof} By replacing $f$ with the function
$(x,y) \to f(x,y)-f(1,y)$ we may assume that $f|(\sigma_1 \cap
\sigma_2) = 0$. (Note that $(x,y) \to f(1,y)$ is always in
$\AC(\sigma)$.)

Suppose first that $f_2 = 0$. Fix $\epsilon > 0$. As $f_1 \in
\AC(\sigma_1)$ there exists $p \in \CTPP(\sigma_1)$ with
$\norm{f_1-p}_{\BV(\sigma_1)} < \epsilon/4$. By the definition of
$\CTPP(\sigma_1)$ there is a triangulation $\{A_i\}_{i=1}^n$ of
$\sigma_1$ such that $p|A_i$ is planar (see \cite[Section 4]{AD}).
Note that $b(y)= p(1,y)$ is a piecewise linear function on $[0,1]$
with $\norm{b}_{\BV[0,1]} = \norm{f_1 - p}_{\BV(\sigma_1 \cap
\sigma_2)} < \epsilon/4$. Extend $p$ to $\sigma_2$ by setting
$p(x,y) = b(y)$. Note that $p \in \CTPP(\sigma)$ and by
\cite[Proposition 4.4]{AD}, $\norm{p|\sigma_2}_{\BV(\sigma_2)} <
\epsilon/4$. Thus, using Theorem~\ref{variation-join},
  \[ \norm{f - p}_{\BV(\sigma)}
     \le \norm{f - p}_{\BV(\sigma_1)} + \norm{f - p}_{\BV(\sigma_2)}
     < \frac{\epsilon}{2}. \]

For arbitrary $f_2$, The same argument will produce a function $q
\in \CTPP(\sigma)$ which approximates to within $\epsilon/2$ the
function which is $f_2$ on $ \sigma_2$ and zero on $\sigma_1$. Thus
the piecewise planar function $p+q$ approximates $f$ to within
$\epsilon$ on $\sigma$. It follows that $f \in \AC(\sigma)$. The
norm estimate is given by Theorem~\ref{variation-join}.
\end{proof}

The conditions on $\sigma_1$ and $\sigma_2$ in
Lemma~\ref{pasting-lemma} could be relaxed considerably. Since we
will not need this greater generality in this paper, we have not
attempted to determine the most general conditions on these sets for
which the above proof works. It is worth noting that one does need
\emph{some} conditions on $\sigma_1$ and $\sigma_2$ or else the
pasted function need not even be of bounded variation.

A major issue in much of this paper will be whether one can always
extend an $\AC(\sigma)$ function to a larger domain.

\begin{quest}\label{extension-quest} Suppose that $\sigma_1 \subseteq \sigma_2$ are
nonempty compact sets. Does there exist $C = C(\sigma_1,\sigma_2)$
such that for every $f \in \AC(\sigma_1)$ there exists ${\tilde f}
\in \AC(\sigma_2)$ such that ${\tilde f}|\sigma_1 = f$ and
  $ \bigl\Vert {\tilde f} \bigr\Vert_{\BV(\sigma_2)} \le C
  \norm{f}_{\BV(\sigma_1)}$?
\end{quest}

The following special case will be needed in Section~\ref{lbl:309} to show that $\AC(\sigma)$ operators are decomposable.

\begin{thm}\label{fill-in-square}
 Let $\sigma$ denote the closed square $[0,1] \times
[0,1]$, and let $\partial \sigma$ denote the boundary of $\sigma$.
Suppose that $b \in \AC(\partial \sigma)$. Then there exists $f \in
\AC(\sigma)$ such that $f|\partial \sigma = b$ and $\normbv{f} \le
28 \norm{b}_{\BV(\partial \sigma)}$.
\end{thm}

\begin{proof} Recall that by \cite[Proposition 4.4]{AD},
if $h \in \AC[0,1]$ is any absolutely continuous function of one
variable, then its extension to the square, $\widehat h (x,y) =
h(x)$, is in $\AC(\sigma)$ with $\ssnorm{\widehat h} =
\norm{h}_{\BV[0,1]}$.

Define $f_s: \sigma \to \mC$ by $f_s(x,y) = (1-y)\,b(x,0)$. Since
$f_s$ is the product of $\AC$ functions of one variable, it is
absolutely continuous on $\sigma$ and
  \[ \normbv{f_s} \le 2 \norm{b(\cdot,0)}_{\BV[0,1]}
                  \le 2 \norm{b}_{\BV(\partial \sigma)}. \]
Similarly, we define
  \begin{align*}
  f_e(x,y) &= (1-x)\, b(0,y),\\
  f_n(x,y) &= y\, b(x,1), \\
  f_w(x,y) &= x\, b(1,y).
  \end{align*}
Let $g = f_s+f_e+f_n+f_w$. Then $g \in \AC(\sigma)$ and $\normbv{g}
\le 8 \norm{b}_{\BV(\partial \sigma)}$.

Let $\Delta_\ell = \{(x,y) \,:\, 0\le y \le x \le 1\}$ and $\Delta_u
= \{(x,y) \,:\, 0 \le x \le y \le 1\}$ denote the lower and upper
closed triangles inside $\sigma$. Now let $p_\ell$ be the affine
function determined by the condition that it agrees with $b-g$ at
the points $(0,0),(1,0)$ and $(1,1)$. Similarly, let $p_u$ be the
affine function which agrees with $b-g$ at the points $(0,0),(0,1)$
and $(1,1)$. Note that $p_\ell(x,x) = p_u(x,x)$ for all $x$. Let
  \[ p(x,y) = \begin{cases}
       p_\ell(x,y),   & (x,y) \in \Delta_\ell, \\
       p_u(x,y),      & (x,y) \in \Delta_u.
       \end{cases}
       \]
Then $p \in \CTPP(\sigma) \subseteq \AC(\sigma)$. Now (using the
facts about $\AC(\sigma)$ functions which only vary in one
direction)
  \[ \var(p,\Delta_\ell) \le \max\{ |p(0,0)-p(1,0)|,
  |p(0,0)-p(1,1)|, |p(1,0)-p(1,1)| \}. \]
Note that
  \begin{align*}
   |p(0,0)-p(1,0)| &\le |b(0,0)-b(1,0)| + |g(0,0)-g(1,0)| \\
         & \le \var(b,\partial \sigma) + \var(g,\sigma) \\
         & \le 9 \norm{b}_{\BV(\partial \sigma)} .
   \end{align*}
This bound also holds for the other terms and hence
$\norm{p}_{\BV(\Delta_\ell)} \le 10 \norm{b}_{\BV(\partial
\sigma)}$. Applying the same argument in the upper triangle, and
then using Theorem~\ref{variation-join} gives that $\normbv{p} \le
20 \norm{b}_{\BV(\partial \sigma)}$.

Let $f = g+p$. Clearly $f \in \AC(\sigma)$ and $\normbv{f} \le 28
\norm{b}_{\BV(\partial \sigma)}$.  Note that $f_e(x,0), f_n(x,0),
f_w(x,0)$ and $p(x,0)$ are all affine functions of $x$, and hence
$f(x,0) - b(x,0)$ is an affine function. But $f(0,0) =
g(0,0)+b(0,0)-g(0,0) = b(0,0)$ and $f(1,0) = b(1,0)$ and so it
follows that $f(x,0) = b(x,0)$ for all $x \in [0,1]$. Similar
arguments hold for the remaining three sides and so $f|\partial
\sigma = b$ as required.
\end{proof}

At the expense of lengthening the reasoning, one could reduce the
constant $28$ in the above theorem. It would be interesting to know
the optimal constant; it seems unlikely that the above construction
would provide this.

In building up $\AC$ functions in Section~\ref{support}, we shall need to make use of the following straightforward extension lemma.

\begin{lem}\label{edges-of-square}
Let $\sigma$ denote the boundary of the square $[0,1]\times[0,1]$.
Denote the four edges of the square as $\{\sigma_i\}_{i=1}^4$. Let
$J$ be a nonempty subset of $\{1,2,3,4\}$ and let $\sigma_J =
\cup_{i \in J} \sigma_i$. Then given any $b \in \AC(\sigma_J)$ there
exists $\hat{b} \in \AC(\sigma)$ with $\hat{b}|\sigma_J = b$ and
$\ssnorm{\hat{b}}_{\BV(\sigma)} \le 4 \norm{b}_{\BV(\sigma_J)}$.
\end{lem}

\begin{proof} Let $T$ denote the circle passing through the 4
vertices of $\sigma$, and let $\pi$ denote the map from $\sigma$ to
$T$ defined by projecting along the rays coming out of the centre of
$\sigma$. Consider a finite list of points $S = \{z_1,\dots,z_n\}
\subseteq \sigma$ with corresponding path $\gamma_S =
\Pi(z_1,\dots,z_n)$. Choose a line $\ell$ in $\mC$ for which
$\gamma_S$ has $\vf(\gamma_S)$ entry points on $\ell$. Note that you
can always do this with $\ell$ passing through the interior of
$\sigma$ and hence $\ell$ is determined by two points $w_1,w_2 \in
\sigma$. Let $\ell_\pi$ denote the line through $\pi(w_1)$ and
$\pi(w_2)$. Since the projection $\pi$ preserves which side of a
line points lie on, $\gamma_{\pi(S)}$ has $\vf(\gamma_S)$ entry
points on $\ell_\pi$.  Conversely, if $\gamma_{\pi(S)}$ has
$\vf(\gamma_{\pi(S)})$ entry points on a line $\ell$, then $\gamma$
must have at least $\vf(\gamma_{\pi(S)})/2$ entry points on the
inverse image of $\ell$ under $\pi$. (The factor of $\frac{1}{2}$
comes from the fact the inverse image of $\ell$ may lie along one of
the edges of $\sigma$.) It follows then that
   \begin{equation}\label{square-to-circle}
    \frac{1}{2}\, \rho(\gamma_S) \le \rho(\gamma_{\pi(S)}) \le
   \rho(\gamma_S).
   \end{equation}

Suppose then that $f \in \BV(\sigma)$. Let $f_\pi: T \to \mC$ be
$f_\pi = f \circ \pi^{-1}$. From (\ref{square-to-circle}) it is
clear that
  \[ \frac{1}{2} \var(f_\pi,T) \le \var(f,\sigma) \le \var(f_\pi,T)
  \]
and so $f_\pi \in \BV(T)$. The same estimate holds when comparing
the variation of $f \in \BV(\sigma_J)$ and that of $f_\pi$ on the
corresponding subset $T_J$ of $T$. But, by \cite[Corollary
5.6]{AD2}, $\BV(T)$ is $2$-isomorphic to the subset of $\BV[0,1]$
consisting of functions which agree at the endpoints. In this final
space, one can extend an $\AC$ function from a finite collection of
subintervals $K$ to the whole of $[0,1]$ by linear interpolation,
without increasing the norm. Note that absolute continuity is
preserved by the isomorphisms between these function spaces. The
factor $4$ comes from collecting together the norms along the
following composition of maps
\[
\begin{CD}
\AC(\sigma_J)     @.       \AC(\sigma) \\
@V{2}V{\pi}V                @A{1}A{\pi^{-1}}A \\
\AC(T_J)          @.       \AC(T) \\
@V{2}VV                      @A{1}AA \\
\AC(K) @>1>\hbox{extend}>  \AC[0,1]
\end{CD}
\]
\end{proof}

Note that if $\sigma_J$ consists of either one side, or else $2$
contiguous sides, then one may extend $b$ to all of $\sigma$ without
increasing of norm using \cite[Proposition~4.4]{AD}. We do not know
whether this is true if, for example,  $\sigma_J$ consists of $2$
opposite sides of the square.

%
%

\section{$\AC(\sigma)$ operators: definition and examples}
\label{lbl:id-s3}

\begin{defn} Suppose that $\sigma \subseteq \mC$ is a nonempty compact set and
that $T$ is a bounded operator on a Banach space $X$. We say that
$T$ is an $AC(\sigma)$ operator if $T$ admits an bounded
$AC(\sigma)$ functional calculus. That is, $T$ is an $\AC(\sigma)$
operator if there exists a bounded unital Banach algebra
homomorphism $\psi: \AC(\sigma) \to B(X)$ for which $\psi(\lam) =
T$.
\end{defn}

Where there seems little room for confusion we shall often say that
$T$ is an $\AC(\sigma)$ operator where one should more properly say
that $T$ is an $\AC(\sigma)$ operator \emph{for some $\sigma$}.

Before proceeding to give some of the general properties of
$AC(\sigma)$ operators, it is appropriate to give the reader some
idea of how this class is related to other standard classes of
operators which arise in spectral theory.

\begin{ex} \label{lbl:189} {\rm
Let $H$ be a Hilbert space and let $T \in B(H)$ be normal. Then $T$
has a $C(\sigma(T))$ functional calculus $\psi$. Then $\psi \vert
\AC(\sigma(T))$ is a linear homomorphism from $\AC(\sigma(T))$ into
$B(X)$. Furthermore $\norm{\psi(f)} \leq \norm{\psi}\norminf{f} \leq
\norm{\psi}\norm{f}_{BV(\sigma(T))}$ for all $f \in \AC(\sigma)$ and
so $\psi \vert \AC(\sigma(T))$ is continuous from $\AC(\sigma(T))$
into $B(H)$. Hence $T$ is an $\AC(\sigma(T))$ operator. Indeed, by
the same argument any scalar type spectral operator (or even
scalar-type prespectral operator) $T$ on a Banach space $X$ is also
an $\AC(\sigma(T))$ operator. (See \cite{hD} for the definitions of
these latter classes of operators.)}
\end{ex}

The operators in the previous example are associated with spectral
expansions which are of an unconditional nature. The motivation for
the present theory is of course to cover operators such as
well-bounded operators, which admit less constrained types of
spectral expansion.

\begin{lem} \label{lbl:909}
Let $T \in B(X)$ be an $\AC(\sigma)$ operator. Suppose that
$\sigma \subset \sigma'$ where $\sigma' \subset \mC$ is compact.
Then $T$ is an $\AC(\sigma')$ operator.
\end{lem}

\begin{proof}
Let $\psi$ be a $\AC(\sigma)$ functional calculus for $T$. Define
$\psi_{\sigma'} : \AC(\sigma') \rightarrow B(X) : f \mapsto \psi(f
\vert \sigma)$. Then $\psi_{\sigma'}$ is a unital linear
homomorphism. Furthermore $\psi_{\sigma'}(\lam) = \psi(\lam \vert
\sigma) = T$. Finally we note from the inequality $\norm{f \vert
\sigma}_{\BV(\sigma)} \leq \norm{f}_{\BV(\sigma')}$ that
$\psi_{\sigma'}$ is continuous. Hence $\psi_{\sigma'}$ is an
$\AC(\sigma')$ functional calculus for $T$.
\end{proof}

The following result was announced in \cite[Section 2]{AD}.

\begin{prop} \label{lbl:508}
Let $T \in B(X)$. The following are equivalent.
\begin{enumerate}
    \item $T$ is well-bounded,
    \item $T$ is an $\AC(\sigma)$ operator for some $\sigma
        \subset \mR$,
    \item $\sigma(T) \subset \mR$ and $T$ is an $\AC(\sigma(T))$
        operator.
\end{enumerate}
\end{prop}

\begin{proof}
Trivially (3) implies (2). Lemma \ref{lbl:909} shows that (2)
implies (1). Say $T$ is well-bounded with functional calculus $\psi
: \AC(J) \rightarrow B(X)$ for some interval $J$. In \cite{AD} we
define a linear isometry $\iota : \AC(\sigma(T)) \rightarrow
\AC(J)$. Define $\psi_{\sigma(T)} : \AC(\sigma(T)) \rightarrow B(X)
: f \mapsto \psi(\iota(f))$. We show that $\psi_{\sigma(T)}$ is an
$\AC(\sigma(T))$ functional calculus for $T$ which will complete the
proof. Clearly $\psi_{\sigma(T)}$ is linear and continuous.
Furthermore, since $\iota(\lam \vert \sigma(T)) = \lam$, we have
that $\psi_{\sigma(T)}(\lam) = T$. To see that $\psi_{\sigma(T)}$ is
a homomorphism we note that if $f, g \in \AC(\sigma(T))$ then
$(\iota(f g) - \iota(f) \iota(g))(\sigma(T)) = \set{0}$. Theorem
4.4.4 of \cite{bA} says we can find a sequence
$\set{h_n}_{n=1}^\infty \subset \AC(J)$ such that $\lim_n \norm{h_n
- (\iota(f g) - \iota(f) \iota(g))}_{\BV(J)} = 0$ and such that for
each $n$, $h_n$ is zero on a neighbourhood of $\sigma(T)$.  This
last condition, by Proposition 3.1.12 of \cite{CF}, implies that
$\psi(h_n) = 0$ for all $n$. Hence $\psi(\iota(f g) - \iota(f)
\iota(g)) = \lim_n \psi(h_n) = 0$, which shows that
$\psi_{\sigma(T)}$ is a homomorphism as claimed.
\end{proof}

As a result of the last proposition we prefer to use the term `real
$\AC(\sigma)$ operator' rather than the term well-bounded operator.
As well as being less descriptive, the term well-bounded operator
also suffers from the fact that it is used for quite a different
concept in the local theory of Banach spaces (see \cite{MTJ} for
example.) We shall however stick with the traditional term for the
remainder of this paper.

The next theorem shows that some important classes of $\AC$
operators are also $\AC(\sigma)$ operators.

\begin{thm} \label{lbl:732}
Let $A \in B(X)$ be well-bounded with functional calculus $\psi :
AC(J) \rightarrow B(X)$ for some interval $J$. Let $f \in AC(J)$
be such that $\rho(f) > 0$. Then $\psi(f)$ is an $\AC(f(J))$
operator.
\end{thm}

\begin{proof}
Define $\psi_f : \AC(f(J)) \rightarrow B(X) : g \mapsto \psi(g \circ
f)$. Then $\psi_f$ is a unital linear homomorphism and $\psi_f(\lam)
= \psi(f)$. By Proposition \ref{lbl:996}, $\psi_f$ is continuous.
\end{proof}

\begin{cor} \label{lbl:744}
Let $A \in B(X)$ be well-bounded and $p$ be a polynomial of one
variable. Then $p(A)$ is an $\AC(p(\sigma(A)))$ operator.
\end{cor}

\begin{cor} \label{lbl:213}
Let $A \in B(X)$ be a well-bounded operator. Then $\exp(i A)$ is
an $\AC(i \exp(\sigma(A)))$ operator.
\end{cor}

We noted earlier that the trigonometrically well-bounded operators
are those operators which can be expressed in the form $\exp(i A)$
where $A \in B(X)$ is a well-bounded operator of type (B). (Indeed one can also insist that
$\sigma(A) \subset [0, 2 \pi]$.) As usual, we denote the unit circle in $\mC$ by $\mT$.

\begin{cor}[\cite{AD2}, Theorem 6.2] \label{lbl:829}
If $T \in B(X)$ be trigonometrically well-bounded then $T$ is an
$\AC(\mT)$ operator. Indeed, if $X$ is reflexive, 
then $T$ is trigonometrically well-bounded operator if and only if it
is an $AC(\mT)$ operator.
\end{cor}

We end this section with a more concrete example.

\begin{ex} {\rm Suppose that $1 < p < \infty$ and that $X$ is the usual
Hardy space $H^p(\mD)$ of analytic functions on the unit disk.
Consider the unbounded operator $Af(z) = z f'(z)$, $f \in H^p(\mD)$
(with natural domain $\{f \,:\, Af \in H^p(\mD)\})$. This operator
arises, for example, as the analytic generator of a semigroup of
composition operators, $T_tf(z) = f(e^{-t} z)$; see \cite{Si}, which
includes a summary of many of the spectral properties of $A$. The
spectrum of $A$ is $\sigma(A) = \mN = \{0,1,2,\dots\}$ with the
corresponding spectral projections $P_k(\sum a_n z^n) = a_k z^k$ ($k
\in \mN$) giving just the usual Fourier components. Suppose then
that $\mu \not\in \sigma(A)$. The resolvent operator $R(\mu,A) =
(\mu I - A)^{-1}$ is a compact operator with spectrum
$\sigma(R(\mu,A)) = \Bigl\{\frac{1}{\mu - k}\Bigr\}_{k=0}^\infty
\cup \{0\}$. From \cite[Theorem 3.3]{CD} it follows easily from the
properties of Fourier series that if $x \in \mR\setminus\mN$, then
$R(x,A)$ is well-bounded. If we fix such an $x$ and take $\mu
\not\in \mR$, then $R(\mu,A) = f(R(x,A))$ where $f(t) =
t/(1+(\mu-x)t)$ is a M{\" o}bius transformation. If $J$ is any
compact interval containing $\sigma(R(x,A))$ then $\rho(f(J)) =
\frac{1}{2}$. Thus $R(\mu,A)$ is an $AC(f(J))$ operator. Thus, all
the resolvents of $A$ are compact $AC(\sigma)$ operators (for some
$\sigma$). Note that none of the resolvents is scalar-type spectral
unless $p=2$.}
\end{ex}

%
%

\section{Properties of $\AC(\sigma)$ operators}
\label{lbl:309}

All $\AC(\sigma)$ operators belong to the larger class of
decomposable operators (in the sense of \cite{CF2}). This will follow immediately from the requirement that the functional calculus map $\psi: \AC(\sigma) \to B(X)$ be what Colojoar{\v a} and Foia{\c s}  term an $\AC(\sigma)$-spectral function. Recall that by Proposition~\ref{lbl:656}, $\AC(\sigma)$ is an admissible algebra.

Suppose that $f \in \AC(\sigma)$. Let $\Omega_f \subseteq \mC$ be the open set $\mC \setminus \supp f$. By Proposition~\ref{lbl:656}, $\Phi_f(\xi) = f_\xi$ is a well-defined map from $\Omega_f$ to $\AC(\sigma)$. 

Following \cite[Section 3.1]{CF2}, the functional calculus map $\psi: \AC(\sigma) \to B(X)$ is an \textbf{$\AC(\sigma)$-spectral function} if, for all $f \in \AC(\sigma)$, the map 
$\psi\circ \Phi_f: \Omega_f \to B(X)$ is analytic on $\Omega_f$.

Since $\psi$ is linear, it suffices to show that the map $\Phi_f$ is differentiable at each point $\xi_0 \in \Omega_f$. To establish this we shall need a technical lemma.

As in \cite{AD3}, let $|x+iy|_\infty = \max(|x|,|y|)$. For $\xi_0 \in \mC$ and $\delta > 0$ let
  \[ B_\infty(\xi_0,\delta) = \{z \in \mC \,:\, |\xi_0 - z|_\infty < \delta\}. \]

\begin{lem}\label{resolvents}
Suppose that $f \in \AC(\sigma)$, $\xi_0 \in \Omega_f$ and that $\delta >0$ is chosen so that $B_\infty(\xi_0,3\delta) \subseteq \Omega_f$. Then there exists a constant $C(\delta,\sigma)$ such that for all $\xi \in B_\infty(\xi_0,\delta)$, there exists $r_\xi \in \AC(\sigma)$ which satisfies
\begin{enumerate}
 \item $r_\xi(z) = \frac{1}{\xi-z}$ for all $z \in \sigma \setminus B_\infty(\xi_0,2\delta)$, and
 \item $\norm{r_\xi}_{\AC(\sigma)} \le C(\delta,\sigma)$.
\end{enumerate}
\end{lem}

\begin{proof}
Suppose first that $\xi_0 \in \sigma$. (The case where $\xi_0 \not\in \sigma$ is similar, but with slightly different norm bounds. The details are left to the reader) 

Let $\sigma_0$ denote the smallest closed square (with sides parallel to the axes) containing $\sigma$ and $B_\infty(\xi_0,3\delta)$.  Let $\sigma_1$ denote the $\overline B_\infty(\xi_0,2\delta)$ and let $\sigma_2$ denote $\sigma_0 \setminus B_\infty(\xi_0,2\delta)$. 

Suppose that $\xi \in B_\infty(\xi_0,\delta)$. The function $z \mapsto \xi - z$ is absolutely continuous on $\sigma_2$ with variation equal to $d = d(\sigma,\delta)$, the length of the diagonal of $\sigma_0$. Since $|\xi - z| \ge \delta$ on $\sigma_2$, Lemma~\ref{lem:ac compact:4590} implies that $r_\xi: \sigma_2 \to \mC$, $z \mapsto (\xi-z)^{-1}$ is in $\AC(\sigma_2)$ with 
  \[ \norm{r_\xi}_{\AC(\sigma_2)} \le \frac{1}{\delta} + \frac{d}{\delta^2}. \]

Clearly $\partial \sigma_1 \subseteq \sigma_2$ so by \cite[Lemmas 3.9 and 4.5]{AD},  $r_\xi|\partial \sigma_1 \in \AC(\partial \sigma_1)$. Using Theorem~\ref{fill-in-square}, we can extend $r_\xi$ to $\sigma_1$ so that $r_\xi| \sigma_1 \in \AC(\sigma_1)$ and $\norm{r_\xi}_{\AC(\sigma_1)} \le 28 \norm{r_\xi}_{\AC(\partial \sigma_1)} \le 28 \norm{r_\xi}_{\AC(\sigma_2)}$. 

By splitting $\sigma_0$ into $9$ smaller rectangles and then repeatedly using Lemma~\ref{pasting-lemma}, one can deduce that $r_\xi \in \AC(\sigma_0)$, and that one has a bound on $\norm{r_\xi}_{\AC(\sigma_0)}$ which depends only on $\sigma$ and $\delta$. Taking the restriction of this function to the original domain $\sigma$ completes the construction.
\end{proof}

\begin{prop}\label{ac-spectral-function}
The functional calculus map $\phi$ for an $\AC(\sigma)$ operator $T \in B(X)$ is an \textbf{$\AC(\sigma)$-spectral function}.
\end{prop}

\begin{proof}
Fix $f \in \AC(\sigma)$, $\xi_0 \in \Omega_f$ and $\delta >0$  so that $B_\infty(\xi_0,3\delta) \subseteq \Omega_f$. Using Lemma~\ref{resolvents}, choose a family of functions $r_\xi$ for $\xi \in B_\infty(\xi_0,\delta)$. Note that $\Phi_f(\xi) = r_\xi f \in \AC(\sigma)$. Thus
\begin{align*}
 \frac{\Phi_f(\xi) - \Phi_f(\xi_0)}{\xi-\xi_0}
     &= \frac{(r_\xi - r_{\xi_0})f}{\xi-\xi_0} \\
     &= -r_\xi r_{\xi_0} f \\
     &= -r_{\xi_0}^2f + r_{\xi_0} (r_{\xi_0}-r_\xi) f \\ 
     &= -r_{\xi_0}^2f + r_{\xi_0} (\xi-\xi_0) r_\xi r_{\xi_0} f \\
     &\to -r_{\xi_0}^2f
\end{align*}
as $\xi \to \xi_0$, by the uniform bound on the norms of the functions $r_\xi$. Composing $\Phi_f$ with the linear map $\psi$ preserves differentiability so $\psi \circ \Phi_f: \Omega_f \to B(X)$ is analytic.
\end{proof}

\begin{prop} \label{lbl:741}
Let $T \in B(X)$ be an $\AC(\sigma)$ operator. Then
\begin{enumerate}
\item $\sigma(T) \subseteq \sigma$.
\item $T$ is decomposable.
\end{enumerate}
\end{prop}

\begin{proof}
This follows from Proposition~\ref{ac-spectral-function} using Theorems 3.1.6 and 3.1.16 in \cite{CF2}.
\end{proof}

In general it is easy to pass between spectral properties of an
operator $T$ and those of affine translations of $T$. One of the
main motivations for developing this theory was to provide a
suitably broad class of operators which is closed under such
transformations. From Theorem \ref{lbl:408} we get the following.

\begin{thm} \label{lbl:814}
Let $T \in B(X)$ be an $\AC(\sigma)$ operator. Let $\alpha, \beta
\in \mC$.  Then $\alpha T + \beta I$ is an $\AC(\alpha \sigma +
\beta)$ operator.
\end{thm}

\begin{proof}
Let $\theta : \AC(\sigma) \rightarrow \AC(\alpha \sigma + \beta)$ be
the isomorphism of Theorem \ref{lbl:408}. Let $\psi$ be the
$\AC(\sigma)$ functional calculus for $T$. Then it is routine to
check that the map $\psi_{\alpha, \beta} : \AC(\alpha \sigma +
\beta) \rightarrow B(X) : f \mapsto \psi(\theta^{-1}(f))$ is an
$\AC(\alpha \sigma + \beta)$ functional calculus for $\alpha T +
\beta I$.
\end{proof}

\begin{thm} \label{lbl:329}
Let $T \in B(X)$ be an $\AC(\sigma)$ operator. Then $T = R + i S$
where $R, S$ are commuting well-bounded operators. Further,
$\sigma(R) = \Real(\sigma(T))$ and $\sigma(S) = \Imag(\sigma(T))$.
\end{thm}

\begin{proof}
Let $\psi$ be an $\AC(\sigma)$ functional calculus for $T$. In
\cite{AD} it is shown in Proposition 5.4 that the map $u :
\AC(\Real(\sigma)) \rightarrow AC(\sigma)$ defined by $u(f)(z) =
f(\Real(z))$ is a norm-decreasing linear homomorphism. Then the map
$\psi_{\Real(\sigma)} : \AC(\Real(\sigma)) \rightarrow B(X) : f
\mapsto \psi(u(f))$ is a continuous linear unital homomorphism.
Hence $R := \psi_{\Real(\sigma)}(\lam \vert \Real(\sigma)) =
\psi(\Real(\lam))$ is well-bounded. Similarly $S :=
\psi(\Imag(\lam))$ is well-bounded. Then $T = \psi(\lam) =
\psi(\Real(\lam) + i\,\Imag(\lam)) := R + i S$. Finally we note that
$R$ and $S$ commute since $\AC(\sigma)$ is a commutative algebra and
$\psi$ is a homomorphism.

The identification of $\sigma(R)$ and $\sigma(S)$ follows
immediately from the spectral mapping theorem  \cite[Theorem 3.2.1]{CF2}
\end{proof}

Splittings which arise from an $\AC(\sigma)$ functional calculus
we call \emph{functional calculus splittings}.

\begin{cor} \label{lbl:781}
The $\AC(\sigma)$ operators are a proper subset of the $\AC$
operators of Berkson and Gillespie.
\end{cor}

\begin{proof}
We note that not all $\AC$ operators are $\AC(\sigma)$ operators.
Example 4.1 of \cite{BDG} shows that the class of $\AC$ operators is
not closed under multiplication by scalars even on Hilbert spaces.
\end{proof}

Not all splittings into commuting real and imaginary well-bounded
parts arise from an $\AC(\sigma)$ functional calculus. This was
shown in the next example which first appeared in \cite{BDG}.

\begin{ex} \label{lbl:690} {\rm
Let $X = L^\infty[0, 1] \oplus L^1[0, 1]$. Define $A \in B(X)$ by $A(f, g) =
(\lam f, \lam g)$. It is not difficult to see that $A$ is well-bounded and
that $\sigma(A) = [0, 1]$. Let $T = (1 + i)A = A + i A$. By Theorem
\ref{lbl:814}, $T$ is an $\AC(\sigma(T))$ operator where $\sigma(T)$ is the
line segment from $0$ to $1 + i$.

The operator $T$ has an infinite number splittings. Define $Q \in
B(X)$ by $Q(f, g) = (0, f)$. In \cite{BDG} it is shown that $A +
\alpha Q$ is well-bounded for any $\alpha \in \mC$. But then $T =
A + i A = A + Q + i(A + i Q)$.

The second splitting cannot come from an $\AC(\sigma)$ functional
calculus. Say $T$ has an $\AC(\sigma)$ functional calculus $\psi$.
Since $\sigma(T)$ is a line segment we can use similar reasoning as
to that in Proposition \ref{lbl:508} to conclude that if $f \in
\AC(\sigma)$ is such that $f(\sigma(T)) = \set{0}$ then $\psi(f) =
0$. Hence if $g \vert \sigma(T) = h \vert \sigma(T)$ then $\psi(g) =
\psi(h)$. In particular since $\Real(\lam) \vert \sigma(T) =
\Imag(\lam) \vert \sigma(T)$ we can only have $\AC(\sigma)$
functional calculus splittings of the form $T = R + i R$. }
\end{ex}

We do not know if it is possible to have several splittings each
arising from an $\AC(\sigma)$ functional calculus. The following
tells us to what extent we can expect splittings to be unique.

\begin{prop} \label{lbl:337}
Let $T \in B(X)$ be an $\AC(\sigma)$ operator. Suppose that $T =
R_1 + i S_1 = R_2 + i S_2$ where $R_1, S_1$ and $R_2, S_2$ are
pairs of commuting well-bounded operators. Then $R_1$ and $R_2$
are quasinilpotent equivalent in the sense of \cite{CF} (as is
$S_1$ and $S_2$). Suppose that $\set{R_1, S_1, R_2, S_2}$ is a
commuting set. Then $(R_1 - R_2)^2 = (S_1 - S_2)^2 = 0$.
Furthermore suppose that $\set{R_1, S_1, R_2, S_2}$ are all type
(B) well-bounded operators. Then $R_1 = R_2$ and $S_1 = S_2$.
\end{prop}

\begin{proof}
This is Theorem 3.2.6 of \cite{CF2} and Theorem 3.7 of \cite{BDG}.
\end{proof}

%
%

\section{The support of the functional calculus}\label{support}

Suppose that $\psi: \AC(\sigma) \to B(X)$ is the functional calculus map for
an $\AC(\sigma)$ operator $T$. The support of $\psi$ is defined as the
smallest closed set $F \subseteq \mC$ such that if 
$\mathrm{supp}f \cap F = \emptyset$, then $\psi(f) = 0$.
It follows from Theorem~\ref{ac-spectral-function} and
\cite[Theorem~3.1.6]{CF2} that the support of $\psi$ is $\sigma(T)$.

It is natural therefore to ask whether such an operator $T$ must
admit an $\AC(\sigma(T))$ functional calculus. By Proposition
\ref{lbl:508}, this is certainly the case if $T$ is well-bounded (that is, if $\sigma(T) \subseteq \mR$),
but the general case remains open.

We shall now give a partial answer to this question, and show that
one may always at least shrink $\sigma$ down to be a compact set not much
bigger than $\sigma(T)$.

\begin{defn} A set $G \subseteq \mC$ is said to be gridlike if it is
a closed polygon with sides parallel to the axes.
\end{defn}

Note that we do not require that a gridlike set be convex, or even
simply connected.

\begin{prop}\label{quotient-prop}
Suppose that $V$ is a gridlike set, that $\sigma$ is compact and
that $V \subseteq \sigma$. Let $I_V = \{f \in \AC(\sigma) \,:\,
\hbox{$f \equiv 0$ on $V$}\}$. Then $\AC(\sigma)/I_V \cong \AC(V)$
as Banach algebras.
\end{prop}

\begin{proof} Define $\Theta: \AC(\sigma)/I_V \to \AC(V)$ by
$\Theta([f]) = f|V$. Then clearly
  \[ \Theta([f]) = \Theta([g]) \iff f|V \equiv g|V \iff f-g \in I_V \]
and so $\Theta$ is well-defined and one-to-one. It is also easy to
see that $\Theta$ is an algebra homomorphism. Since
  \begin{align*}
  \norm{\Theta([f])}
    & = \norm{f|V}_{\BV(V)} \\
    & = \inf_{g \in I_V} \norm{f+g|V}_{\BV(V)} \\
    & \le \inf_{g \in I_V} \normbv{f+g} \\
    &= \norm{[f]}_{\AC(\sigma)/I_V}
  \end{align*}
the map $\Theta$ is bounded.

The hard part of the proof is to show that $\Theta$ is onto. That
is, given $f \in \AC(V)$, there exists $F \in \AC(\sigma)$ so that
$F|V = f$.

Choose then a square $J \times K$ containing $\sigma$. Extending the
edges of $V$ produces a grid on $J \times K$, determining $N$ closed
subrectangles $\{\sigma_k\}_{k=1}^{N}$.

Suppose now that $f \in \AC(V)$. Our aim is to define $\hat{f} \in
\AC(J \times K)$ with $\hat{f}|V = f$ and $\ssnorm{\hat{f}}_{\BV(J
\times K)} \le C \norm{f}_{\BV(V)}$.

Fix an ordering of the rectangles $\{\sigma_k\}$ so that
\begin{enumerate}
\item there exists $k_0$ such that $\sigma_k \subseteq V$ if and
only if $k \le k_0$, and
\item for all $\ell$, $\sigma_\ell$ intersects $\cup_{k < \ell}
\sigma_k$ on at least one edge of $\sigma_\ell$.
\end{enumerate}
Let $E_0$ denote the union of the edges of the rectangles $\sigma_k$
for $k \le k_0$ and let $b$ be the restriction of $f$ to $E_0$. Note
that $b$ is absolutely continuous  on $E_0$ and if $e$ is any edge of
any rectangle $\sigma_k$ ($k \le k_0$), then $b|e \in \AC(e)$ with
$\norm{b|e}_{\BV(e)} \le \norm{b}_{\BV(E_0)} \le
\norm{f}_{\BV(V)}$. Now apply Lemma~\ref{edges-of-square}
to recursively extend $b$ to the set $E$ of all edges of rectangles
$\sigma_k$, $1 \le k \le N$, so that $b \in \AC(E)$ and
$\norm{b}_{\BV(E)} \le C_N \norm{f}_{\BV(V)}$.

For $1 \le k \le k_0$, let $f_k = f|\sigma_k$, so that $f_k \in
\AC(\sigma_k)$ and $\norm{f_k}_{\BV(\sigma_k)} \le
\norm{f}_{\BV(V)}$. Suppose alternatively that $k_0 < k \le N$. By
Theorem~\ref{fill-in-square} we can find $f_k \in \AC(\sigma_k)$
with $f_k|\partial \sigma_k = b|\partial\sigma_k$ and
$\norm{f_k}_{\BV(\sigma_k)} \le 28
\norm{b|\partial\sigma_k}_{\BV(\partial\sigma_k)} \le 28 C_N
\norm{f}_{\BV(V)}$.

Define $\hat{f}:J \times K \to \mC$ such that $\hat{f}|\sigma_k =
f_k$. That $\hat{f}$ is in $\AC(J \times K)$ with
$\ssnorm{\hat{f}}_{\BV(J \times K)} \le 28 C_N N \norm{f}_{\BV(V)}$
follows from Lemma~\ref{pasting-lemma} (first patching together all
the rectangles in each row, and then all the rows together). We can now
let $F = \hat{f}|\sigma$.

It follows then that $\Theta$ is onto and hence is a Banach algebra
isomorphism.
\end{proof}

\begin{thm} \label{lbl:662}
Let $T \in B(X)$ be an $\AC(\sigma)$ operator for some $\sigma
\subset \mC$. Let $U$ be an open neighbourhood of $\sigma(T)$. Then
$T$ is an $\AC(\overline{U})$ operator.
\end{thm}

\begin{proof} Suppose that $T$, $\sigma$ and $U$ are as stated.
Choose a square $J \times K$ containing $U \cup \sigma$. By
Lemma~\ref{lbl:909}, $T$ admits an $\AC(J \times K)$ functional
calculus $\psi$.

Consider an equispaced grid on $J \times K$, determining $n^2$
subsquares $\{\sigma_k\}_{k=1}^{n^2}$. Let $V = V(n)$ be the union
of all those $\sigma_k$ which intersect $\sigma(T)$. For $n$ large
enough
   \[ \sigma(T) \subseteq \mathrm{int}(V) \subseteq V \subseteq  U. \]
For the rest of the proof, fix such an $n$.

As in Proposition~\ref{quotient-prop}, let $I_V = \{ f \in \AC(J
\times K) \,:\, f|V \equiv 0\}$, so that $\AC(J\times K)/I_V \cong
\AC(V)$ via the isomorphism $\Theta$. Note that $I_V \subseteq
\mathrm{ker} (\psi)$ since if $f \in I_V$, then $\mathrm{supp}f \cap
\sigma(T) = \emptyset$. Thus the map $\tilde{\psi}: \AC(J\times
K)/I_V \to B(X)$,
  \[ \tilde{\psi}([f]) = \psi(f) \]
is a well-defined algebra homomorphism with $\ssnorm{\tilde{\psi}}
\le \norm{\psi}$.

We may therefore define $\hat{\psi}: \AC(\overline{U}) \to B(X)$ by
$\hat{\psi}(f) = \tilde{\psi}([\Theta^{-1}(f|V)])$. Note that
$\hat{\psi}$ is a bounded algebra homomorphism and that, since
$\Theta([\lam]) = \lam|V$, $\hat{\psi}(\lam) = \psi(\lam) = T$. Thus
$\hat{\psi}$ is an $\AC(\overline{U})$ functional calculus for $T$.
\end{proof}

\begin{cor}\label{lbl:cor-to-662}
Let $T \in B(X)$ be an $\AC(\sigma_0)$ operator for some compact set
$\sigma_0$. Then
  \[ \sigma(T) = \bigcap \{ \sigma \,:\,
               \hbox{$T$ has an $\AC(\sigma)$ functional
               calculus}\}. \]
\end{cor}

The proof of Theorem~\ref{lbl:662} depends on two vital facts. The
first is that the map $\Theta$ is an isomorphism. The second is that
$I_V \subseteq \mathrm{ker}(\psi)$. To show that every $\AC(\sigma)$
operator is an $\AC(\sigma(T))$ operator, it would suffice to show
that
\begin{enumerate}
  \item\label{Q1} the restriction map $\AC(\sigma) \to \AC(\sigma(T))$, $f
  \mapsto f|\sigma(T)$ is onto. This is basically equivalent to
  answering Question~\ref{extension-quest}.
  \item\label{Q2} given any $f \in \AC(\sigma)$ with $f|\sigma(T) \equiv 0$,
  there exists a sequence $\{f_n\} \subseteq \AC(\sigma)$ with
  $\normbv{f-f_n} \to 0$ and $\mathrm{supp} f_n \cap \sigma(T) =
  \emptyset$ for all $n$.
\end{enumerate}
Proving (\ref{Q1}) and (\ref{Q2}) when $\sigma(T)$ is a complicated
compact set would appear to require new ways of estimating the
two-dimensional variation used in our definitions.

If $T \in B(X)$ is an $\AC(\sigma(T))$ operator then $T$ has
spectral theorems similar to those for normal operators. Recall
from \cite{nD} the definition of the local spectrum $\sigma_T(x)$
of $x \in X$ for an operator $T \in B(X)$ with the single-valued
extension property. From \cite{LV} if $T \in B(X)$ is an
$\AC(\sigma)$ operator (and hence decomposable) then those $x \in
X$ such that $\sigma_T(x) = \sigma(T)$ are second countable in
$X$.

\begin{thm} \label{lbl:091}
Suppose that $T \in B(X)$ is an $\AC(\sigma(T))$ operator with
functional calculus $\psi : \AC(\sigma(T)) \rightarrow B(X)$. Then
$\psi$ is injective. Hence we can identify $\AC(\sigma(T))$ with a
subalgebra of $B(X)$. Furthermore suppose that $x \in X$ is such
that $\sigma_T(x) = \sigma(T)$. Then the map $\AC(\sigma(T))
\rightarrow X : f \mapsto \psi(f)x$ is injective, and so we can
identify $\AC(\sigma(T))$ with a subspace of $X$.
\end{thm}

\begin{proof}
Let $x \in X$ be such that $\sigma_T(x) = \sigma(T)$. To prove the
theorem it suffices to show that if $f \in \AC(\sigma(T))$ and $f
\neq 0$ then $\psi(f)x \neq 0$. Let $\lambda_0 \in \sigma(T)$ be
such that $f(\lambda_0) \neq 0$. Since $f$ is continuous we can find
an open neighbourhood $V$ of $\lambda_0$ such that $0 \not \in
f(V)$. We can choose $g \in \AC(\sigma(T))$ such that $(f g)(V) =
\set{1}$. If we show $\psi(f g) x \neq 0$ this will imply, since
$\psi$ is a homomorphism, that $\psi(f)x \neq 0$. Hence we can
assume that $f(V) = \set{1}$. Let $U$ be an open set such that
$\set{U, V}$ is an open cover of $\sigma(T)$ and such that
$\lambda_0 \not \in U$. By Lemma 5.2.3 of \cite{bA} we can find
non-zero $x_U, x_V \in X$ such that $x = x_U + x_V$ and where
$\sigma_T(x_U) \subset U$ and $\sigma_T(x_V) \subset V$. Since
$\sigma_T(x) \subset \sigma_T(x_U) \cup \sigma_T(x_V)$ we have that
$\lambda_0 \in \sigma_T(x_V)$ and $\lambda_0 \not \in
\sigma_T(x_U)$. Assume that $\psi(f)x = 0$. Then 
$0 = \psi(f)(x_U +x_V) = \psi(f)x_U + x_V$ 
since $f$ is one on $V$. It follows that
$\sigma_T(x_V) = \sigma_T(-\psi(f)x_U) = \sigma_T(\psi(-f)x_U)
\subset \sigma_T(x_U)$. Then we have the contradiction that
$\lambda_0 \in \sigma_T(x_V) \subset \sigma_T(x_U) \not \ni
\lambda_0$. Hence $\psi(f)x \neq 0$.
\end{proof}

Since every $\AC(\sigma)$ operator is also an $\AC$ operator, the
results of \cite{DW} give a representation theorem for compact
$\AC(\sigma)$ operators. Specifically, if $T \in B(X)$ is a compact
$\AC(\sigma)$ operator with nonzero eigenvalues $\{\mu_j\}$ and
corresponding Riesz projections $\{P_j\}$, then
  \begin{equation}\label{comp-sum}
   T = \sum_j \mu_j P_j
  \end{equation}
where the sum converges in norm under a particular specified
ordering of the eigenvalues. Given a sequence of real numbers
$\{\mu_j\}$ and disjoint projections $\{P_j\} \subseteq B(X)$,
necessary and sufficient conditions are known which ensure that the
operator defined via (\ref{comp-sum}) is well-bounded (\cite[Theorem
3.3]{CD}). At present an analogous result for compact $\AC(\sigma)$
operators in unknown. These questions are pursued more fully in
\cite{AD3} where, for example, various sufficient conditions for
(\ref{comp-sum}) to define a compact $\AC(\sigma)$ operator are
given.

%
%

\section{Spectral resolutions}
\label{lbl:SpecRes}

The theory of well-bounded operators is at its most powerful if one
adds the additional assumption that the functional calculus map for
$T$ is `weakly compact'. That is, for all $x \in X$, the map
$\psi_x: \AC(\sigma(T)) \to X$, $f \mapsto \psi(f)x$ is weakly
compact. In this case $T$ admits an integral representation with
respect to a spectral family of projections $\{E(\mu)\}_{\mu \in
\mR}$. The integration theory for spectral families allows one to
define
  \[ f(T) = \widehat{\psi}(f)
          = \int_{\sigma(T)}^\oplus f(\mu) \, dE(\mu) \]
for all $f \in \BV(\sigma)$ giving an extended functional calculus
map. (This integral is more usually written as $\int_{J}^\oplus f(\mu)
\, dE(\mu)$, where $J$ is some compact interval containing
$\sigma(T)$. We have written it in the above form to stress that the
value of the integral only depends on the values of $f$ on
$\sigma(T)$.) If $\psi$ is not weakly compact, then there may be no
spectral resolution consisting of projections on $X$. A suitable
family of projections on $X^*$, known as a decomposition of the
identity, does always exist, but the theory here is much less
satisfactory.

Obviously extending this theory to cover general $\AC(\sigma)$
operators with a weakly compact functional calculus is highly
desirable. At present a full analogue of the well-bounded theory has
not been found, but we are able to show that each such operator does
admit a nice spectral resolution from which the operator may be
recovered. The following definition extends the definition for
well-bounded operators.

\begin{defn} Let $T \in B(X)$ be an $\AC(\sigma)$ operator with
functional calculus map $\psi$. Then $T$ is said to be of type~(B)
if for all $x \in X$, the map $\psi_x: \AC(\sigma(T)) \to X$, $f
\mapsto \psi(f)x$ is weakly compact.
\end{defn}

Obviously every $\AC(\sigma)$ operator on a reflexive Banach space
is of type~(B),
 as is every scalar-type spectral operator on a
general Banach space (see \cite{K}). The weak compactness of the
functional calculus removes one of the potential complications with
studying $\AC(\sigma)$ operators.

\begin{lem} \label{lbl:129}
Let $T \in B(X)$ have a weakly compact $\AC(\sigma)$ functional
calculus. Then it has a unique splitting $T = R + i S$ where $R$ and
$S$ are commuting type (B) well-bounded operators.
\end{lem}

\begin{proof}
Recall if we set $R = \psi(\Real(\lam))$ and $S = \psi(\Imag(\lam))$
then $R$ and $S$ are commuting well-bounded operators. The
$\AC(\Real(\sigma(T)))$ functional calculus for $R$ is given by $f
\mapsto \psi(u(f))$ and so is clearly weakly compact. Hence $R$ is type
(B). Similarly $S$ is type (B). Uniqueness follows from Proposition
\ref{lbl:337}.
\end{proof}

If $T$ is a well-bounded operator of type~(B) with spectral family
$\{E(\mu)\}_{\mu \in \mR}$, then, for each $\mu$, $E(\mu)$ is the
spectral projections for the interval $(-\infty,\mu]$. The natural
analogue of this in the $\AC(\sigma)$ operator setting is to index
the spectral resolution by half-planes. Modelling the plane as
$\mR^2$, each closed half-plane is specified by a unit vector
$\theta \in \mT$ and a real number $\mu$:
  \[ H(\theta,\mu) = \{ z \in \mR^2 \,:\, z\cdot \theta \le \mu \} .\]
Let $\HP$ denote the set of all half-planes in $\mR^2$. The
following provisional definition contains the minimal conditions one
would require of a spectral resolution for an $\AC(\sigma)$
operator.

\begin{defn}\label{HPSF} Let $X$ be a Banach space. A half-plane spectral family on
$X$ is a family of projections $\{E(H)\}_{H \in \HP}$ satisfying:
\begin{enumerate}
  \item\label{HPSF-1} $E(H_1)\, E(H_2) = E(H_2)\, E(H_1)$ for all $H_1,H_2 \in
  \HP$;
  \item there exists $K$ such that $\norm{E(H)} \le K$ for all $H
  \in \HP$;
  \item for all $\theta \in \mT$,
  $\{E(H(\theta,\mu))\}_{\mu \in \mR}$ forms a spectral family of
  projections.
  \item\label{HPSF-4} for all $\theta \in \mT$, if $\mu_1 < \mu_2$, then
  $E(H(\theta,\mu_1))\, E(H(-\theta,-\mu_2)) = 0$.
\end{enumerate}
The radius of $\{E(H)\}$ is the (possibly infinite) value
  \[ r(\{E(H)\}) =  \inf \{ r \,:\,
     \hbox{for all $\theta$, $E(H(\theta,\mu)) = I$ for all $\mu > r$} \} . \]
\end{defn}

Suppose that $\sigma \subset \mR^2$ is a nonempty compact set. Given
any unit direction vector $\theta$, let $\sigma_\theta = \{ z\cdot
\theta \,:\, z \in \sigma\} \subseteq \mR$. Define the subalgebra of
all $\AC(\sigma)$ functions which only depend on the component of
the argument in the direction $\theta$,
  \[ \AC_\theta(\sigma) = \{ f \in \AC(\sigma) \,:\,
  \hbox{there exists $u \in \AC(\sigma_\theta)$ such that $f(z) =
  u(z \cdot \theta)$}\} .\]
By Proposition~3.9 and Lemma~3.10 of \cite{AD}, there is a norm 1
isomorphism $U_\theta: \AC(\sigma_\theta) \to \AC_\theta(\sigma)$.

Let $T \in B(X)$ be an $\AC(\sigma)$ operator of type~(B), with
functional calculus map $\psi$. The algebra homomorphism
$\psi_\theta: \AC(\sigma_\theta) \to B(X)$, $u \mapsto \psi(U_\theta
u)$ is clearly bounded and weakly compact. It follows then from the
spectral theorem for well-bounded operators of type~(B) (see, for
example, \cite{qBD}) that there exists a spectral family
$\{E(H(\theta,\mu))\}_{\mu \in \mR}$, with $\norm{E(H(\theta,\mu))}
\le 2 \norm{\psi}$ for all $\mu$. We have thus constructed a
uniformly bounded family of projections $\{E(H)\}_{H \in \HP}$. To
show that this family is a half-plane spectral family it only
remains to verify (\ref{HPSF-1}) and (\ref{HPSF-4}).

Suppose then that $E_1 = E(\theta_1,\mu_1)$ and $E_2 =
E(\theta_2,\mu_2)$. For $\mu \in \mR$ and $\delta > 0$, let
$g_{\mu,\delta}: \mR \to \mR$ be the function which is $1$ on
$(-\infty,\mu]$, is $0$ on $[\mu+\delta,\infty)$ and which is linear
on $[\mu,\mu+\delta]$. Let $h_\delta = U_{\theta_1}
(g_{\mu_1,\delta})$ and $k_\delta = U_{\theta_2}(
g_{\mu_2,\delta})$. The proof of the spectral theorem for
well-bounded operators shows that $E_1 = \lim_{\delta \to 0^+}
\psi(h_\delta)$ and $E_2 = \lim_{\delta \to 0^+} \psi(k_\delta)$,
where the limits are taken in the weak operator topology in $B(X)$.
Thus, if $x \in X$ and $x^* \in X^*$,
  \begin{align*}
  \ipr<E_1 E_2 x,x^*>
    &= \lim_{\delta \to 0^+} \ipr<\psi(h_\delta) E_2 x,x^*> \\
    &= \lim_{\delta \to 0^+} \ipr< E_2 x, \psi(h_\delta)^*x^*> \\
    &= \lim_{\delta \to 0^+} \left( \lim_{\beta \to 0^+}
         \ipr<\psi(h_\delta) \psi(k_\beta) x,x^*> \right)\\
    &=\lim_{\delta \to 0^+} \left( \lim_{\beta \to 0^+}
         \ipr<\psi(k_\beta) \psi(h_\delta) x,x^*> \right)\\
    &=\lim_{\delta \to 0^+} \ipr<\psi(h_\delta) x, E_2^*x^*> \\
    &= \ipr<E_2 E_1 x,x^*>
  \end{align*}
Verifying (\ref{HPSF-4}) is similar. Fix $\theta \in \mT$ and $\mu_1
< \mu_2$. Let $E_1 = E(\theta,\mu_1)$ and $E_2 = E(-\theta,-\mu_2)$.
Let $h_\delta = U_{\theta} (g_{\mu_1,\delta})$ and $k_\delta =
U_{-\theta}( g_{-\mu_2,\delta})$ so that $E_1 = \lim_{\delta \to
0^+} \psi(h_\delta)$ and $E_2 = \lim_{\delta \to 0^+}
\psi(k_\delta)$. The result follows by noting that for $\delta$
small enough, $h_\delta k_\delta = 0$.

We have shown then that $\{E(H)\}_{H \in \HP}$ is a half-plane
spectral family.

\smallskip
For notational convenience, we shall identify the direction vector
$\theta \in \mR^2$ with the corresponding complex number on the unit
circle. Thus, for example, we identify $(0,1)$ with $i$.

For $\theta \in \mT$, the spectral family $\{E(\theta,\mu)\}_{\mu
\in \mR}$ defines a well-bounded operator of type~(B)
  \begin{equation}\label{T-theta}
   T_\theta = \int_{\sigma_\theta} \mu \, dE(\theta,\mu).
   \end{equation}
Clearly the map $\lam_\theta = z \cdot \theta$ lies in $\AC_\theta(\sigma) \subseteq \AC(\sigma)$ and the construction of the spectral family ensures that $\psi(\lam_\theta) = T_\theta$.
Since $\lam = \theta \lam_\theta + i\theta \lam_{i \theta}$ we have that 
  \begin{equation}\label{theta-decomp}
  T = \theta T_\theta + \,i\theta\, T_{i\theta}.
  \end{equation}
In particular, using Theorem~\ref{lbl:329} and Theorem~\ref{lbl:337} we have
that $T$ has the unique splitting into real and imaginary parts
  \begin{equation}\label{recon-formula}
  T = T_1
        + \,i\, T_i .
  \end{equation}
One consequence of these identities is that $T$ may be recovered from the half-plane spectral family
produced by the above construction.

Note that Theorem~\ref{lbl:329} and the fact that $\theta^{-1} T = T_\theta + i T_{i\theta}$ imply that $\sigma(T_\theta) = \Real(\sigma(\theta^{-1}T))$. Thus, if $r(\cdot)$
denotes the spectral radius, then $r(T_\theta) \le r(\theta^{-1}T) = r(T)$. 

Since there exists $\theta \in \mT$ for
which $r(T_\theta) = r(T)$, we have the following result.

\begin{prop} With $T$ and $\{E(H)\}$ as above, $r(\{E(H)\}) = r(T)$.
\end{prop}

Note that if we define $f_\theta \in \AC(\sigma)$ by $f_\theta(z) =
z\cdot \theta$, then $T_\theta = \psi(f_\theta) = f_\theta(T)$. In
particular, if $\omega = (1/\sqrt{2},1/\sqrt{2})$, then $f_\omega =
(f_1 + f_i)/\sqrt{2}$, and hence
  \[ T_\omega = \psi(f_\omega) = (T_1+T_i)/\sqrt{2}. \]
This proves the following proposition. Note that in general the sum
of two commuting well-bounded operators need not commute.

\begin{prop} Let $T$ be an $\AC(\sigma)$ operator of type~(B), with
unique splitting $T = R+iS$. Then $R+S$ is also well-bounded.
\end{prop}

\begin{quest}
Suppose that $R$ and $S$ are commuting well-bounded operators whose
sum is well-bounded. Is $R+iS$ an $\AC(\sigma)$ operator?
\end{quest}

It is clear that given any half-plane spectral family $\{E(H)\}_{H
\in \HP}$ with finite radius, Equation~(\ref{recon-formula}) defines
$T \in B(X)$ which is an $\AC$~operator in the sense of Berkson and
Gillespie. It is not clear however, that $T$ need be an
$\AC(\sigma)$ operator. In particular, if we define $T_\theta$ via
Equation~(\ref{T-theta}), then it is not known whether the identity
(\ref{theta-decomp}) holds.

\begin{quest}
Is there a one-to-one correspondence between $\AC(\sigma)$ operators
of type~(B) and half-plane spectral families with finite radius? If
not, can one refine Definition~\ref{HPSF} so that such a
correspondence exists?
\end{quest}

%
%

\section{Extending the functional calculus}
\label{lbl:454}

Given a $\AC(\sigma)$ operators of type~(B) its associated
half-plane spectral family (as constructed above), it is natural to
ask whether one can develop an integration theory which would enable
the functional calculus to be extended to a larger algebra than
$\AC(\sigma)$.

The spectral family associated to a well-bounded operator $T$ of
type~(B) allows one to associate a bounded projection with any set
of the form $\bigcup_{j=1}^n \sigma(T) \cap I_j$, where
$I_1,\dots,I_n$ are disjoint intervals of $\mR$. Let $\sigma =
\sigma(T) \subset \mR$ and let $\NP(\sigma(T))$ denote the algebra
of all such sets. It is easy to check the following.

\begin{thm} \label{lbl:777}
Let $T \in B(X)$ be a type (B) well-bounded operator with functional
calculus $\psi$. Then there is a unique map $E : \NP(\sigma(T))
\rightarrow \Proj(X)$ satisfying the following
\begin{enumerate}
    \item $E(\emptyset) = 0$, $E(\sigma(T)) = I$,
    \item $E(A \cap B) = E(A)E(B) = E(B)E(A)$ for all $A, B \in \NP(\sigma(T))$,
    \item $E(A \cup B) = E(A) + E(B) - E(A \cap B)$
        for all $A, B \in \NP(\sigma(T))$,
    \item\label{NP-norm-bound}
    $\normbv{E(A)} \leq \norm{\psi} \normbvj{\chi_A}$
        for all $A \in \NP(\sigma(T))$,
    \item if $S \in B(X)$ is such that $T S = S T$ then $E(A) S = S E(A)$
        for all $A \in \NP(\sigma(T))$,
    \item $\Range(E(A)) = \set{x \in X : \sigma_T(x) \subseteq A}$.
\end{enumerate}
\end{thm}

For general $\AC(\sigma)$ operators, the natural algebra of sets is
that generated by the closed half-planes. This algebra has been
studied in various guises, particularly in the setting of
computational geometry. The sets that can be obtained by starting
with closed half-planes and applying a finite number of unions,
intersections and set complements are sometimes known as Nef
polygons. The set of Nef polygons in the plane, $\NP$, clearly
contains all polygons, lines and points in the plane. For more
information about Nef polygons, or more generally their
$n$-dimensional analogues, Nef polyhedra, we refer the reader to
\cite{BIC}, \cite{HKM} or \cite{Nef}.

Let $\sigma$ be a nonempty compact subset of $\mC$. Define
  \[ \NP(\sigma) = \{ A \,:\, \hbox{$A = \sigma \cap P$ for some $P
  \in \NP$}\}. \]
It is clear that given an $AC(\sigma)$ operator of type~(B), one may
use the half-plane spectral family constructed in the previous
section to associate a projection $E(A) \in B(X)$ with each set $A
\in \NP(\sigma)$. The major obstacle in developing a suitable
integration theory in this setting is in providing an analogue of
condition (\ref{NP-norm-bound}) in Theorem~\ref{lbl:777}.

Note that if $A \in \NP(\sigma)$, then $\chi_A \in \BV(\sigma)$.
Rather than forming $E(A)$ by a finite combination of algebra
operations, one might try to define $E(A)$ directly as we did when
$A$ was a half-plane. That is, one may try to write
  \[ E(A) = WOT-\lim_\alpha \psi(h_\alpha) \]
where $\{h_\alpha\}$ is a suitable uniformly bounded net of
functions in $\AC(\sigma)$  which approximates $\chi_A$ pointwise.
It is shown in \cite{bA} that if $A$ is a closed polygon then this
may be done but only under the bound
  $ \norm{h_\alpha} \le V_A$ .
Here $V_A$ is a constant depending on $A$. This allows one to prove
a weaker version of Theorem~\ref{lbl:777}, with condition
(\ref{NP-norm-bound}) replaced by $\norm{E(A)} \le V_A \norm{\psi}$.
It remains an open question as whether one can do this with  $V_A
\le 2 \norm{\chi_A}$. However, if $A$ is a closed convex polygon
contained in the interior of $\sigma$, then this is possible.

\begin{quest} Does every $\AC(\sigma)$ operator of type~(B) admit a
$\BV(\sigma)$ functional calculus?
\end{quest}

It might be noted in this regard that all the examples of
$\AC(\sigma)$ operators of type~(B) given in Section~\ref{lbl:id-s3}
do admit such a functional calculus extension.

%
%

\end{document}